\documentclass[11pt, letter]{amsart}

\usepackage{graphicx,amscd}
\newtheorem{theorem}{Theorem}[section]

\newtheorem{lemma}[theorem]{Lemma}

\newtheorem*{claim}{Claim}

\theoremstyle{remark}

\def\T{{\mathcal T}}
\def\S{{\mathcal S}}
\def\U{{\mathcal U}}
\def\zed{{\mathbb Z}}

\begin{document}

\title{Examples of reducible and finite Dehn fillings}

\author[S. Kang]{Sungmo Kang}
\address{Department of Mathematics, University of Texas at Austin, 1 University Station C1200 Austin, TX
78712-0257} \email{skang@math.utexas.edu}

\begin{abstract}
If a hyperbolic $3$-manifold $M$ admits a reducible and a finite
Dehn filling, the distance between the filling slopes is known to
be 1. This has been proved recently by Boyer, Gordon and Zhang.
The first example of a manifold with two such fillings was given
by Boyer and Zhang. In this paper, we give examples of hyperbolic
manifolds admitting a reducible Dehn filling and a finite Dehn
filling of every type: \textit{cyclic, dihedral, tetrahedral,
octahedral and icosahedral}.
\end{abstract}

\maketitle

\section{Introduction}
Let $M$ be a compact connected orientable $3$-manifold with a
torus boundary component $\partial_0 M$ and $r$ a {\em slope}, the
isotopy class of an essential simple closed curve, on $\partial_0
M$. The manifold obtained by {\em $r$-Dehn filling} on $M$ is
defined to be $M(r)=M\cup V$, where $V$ is a solid torus glued to
$M$ along $\partial_0 M$ so that $r$ bounds a disk in $V$.

We say that $M$ is \textit{hyperbolic} if $M$ with its torus
boundary components removed has a complete hyperbolic structure of
finite volume. For a pair of slopes $r_1$ and $r_2$ on $\partial_0
M$, the distance $\Delta (r_1,r_2)$ denotes their minimal
geometric intersection number. For a hyperbolic manifold $M$, if
both $M(r_1)$ and $M(r_2)$ fail to be hyperbolic, then the upper
bounds for $\Delta(r_1,r_2)$ have been established in various
cases. See \cite{G3}

We are interested in the case of reducible and finite Dehn
fillings (i.e. a filling whose fundamental group is finite). The
first step towards determining the least upper bound for
$\Delta(r_1,r_2)$ in that case was achieved in \cite{BCSZ}, where
the bound of 2 was established. Furthermore, in \cite{BCSZ} it was
shown that if $\Delta(r_1,r_2)=2$ then $H_1(M) \cong \zed \oplus
\zed_2$ and the reducible Dehn filling is homeomorphic to $L(2,1)
\# L(3,1)$. Recently, this special case was eliminated by Boyer,
Gordon and Zhang \cite{BGZ}, showing that $\Delta(r_1,r_2)=1$

There are five types of finite group that can occur as the
fundamental group of a 3-manifold: \textit{cyclic},
\textit{dihedral}, \textit{tetrahedral}, \textit{octahedral} and
\textit{icosahedral}. See \cite[Section 1]{BZ1} for the definition
of these. It is known that all these types of finite group can be
realized as the fundamental group of either a lens space or a
Seifert fibered space over $S^2$ with three exceptional fibers of
orders $a,b,c$ satisfying $1/a+1/b+1/c>1$. The latter is denoted
by $S^2(a,b,c)$. More precisely cyclic, dihedral, tetrahedral,
octahedral and icosahedral types are the fundamental groups of a
lens space, $S^2(2,2,n), S^2(2,3,3), S^2(2,3,4)$ and $S^2(2,3,5)$
respectively.

The following are the known examples of hyperbolic manifolds
admitting a reducible and a finite Dehn filling so far;

(1) In \cite[Example 7.8]{BZ2}, Boyer and Zhang gave the first
example of a manifold realizing the upper bound 1. This manifold
admits a reducible Dehn filling, a Dehn filling of cyclic type and
a Dehn filling of dihedral type. See Sections 2, 3 for details.

(2) In \cite[Section 4]{EW}, Eudave-Mu\~{n}oz and Wu gave an
infinite family of hyperbolic manifolds admitting a reducible Dehn
filling and a Dehn filling of cyclic type. See Section 2 for
details.

(3) In \cite[Section 4]{Lee2}, Lee constructed a family of
hyperbolic manifolds and showed that these manifolds admit an
$S^1\times S^2$ Dehn filling and a toroidal Dehn filling with
distance 2. However we observe that some manifolds in this family
also admit a Dehn filling of icosahedral type. See Section 6 for
details.

In this paper we give examples of all types of finite group. The
technique of constructing examples is to use tangles and double
branched covers. The double branched cover of a tangle is a
3-manifold whose boundary consists of tori. Also performing a
rational tangle filling on a given tangle corresponds to
performing a Dehn filling to the corresponding double branched
cover since the double branched cover of a rational tangle is a
solid torus. By using tangle arguments, we establish the following
theorem.

\begin{theorem}\label{main}
There are hyperbolic manifolds admitting a reducible Dehn filling
and a finite Dehn filling of every type at distance 1: cyclic,
dihedral, tetrahedral, octahedral and icosahedral.
\end{theorem}

This paper is organized as follows. In Section 2 we construct
another infinite family of hyperbolic manifolds admitting a
reducible Dehn filling (i.e. the connected sum of two lens spaces
of arbitrary orders) and a Dehn filling of cyclic type (i.e. a
lens space). In Section 3 we present hyperbolic manifolds which
admit a reducible Dehn filling (i.e. the connected sum of two lens
spaces) and a Dehn filling of dihedral type (i.e. $S^2(2,2,n)$).
In Section 4 we give the first example of a hyperbolic manifold
admitting a reducible Dehn filling (i.e. the connected sum of two
lens spaces) and a Dehn filling of tetrahedral type (i.e.
$S^2(2,3,3)$). In Section 5 we describe the first examples
(infinitely many) of hyperbolic manifolds which admit a reducible
Dehn filling (i.e. the connected sum of a lens space and a small
Seifert fibered space) and a Dehn filling of octahedral type (i.e.
$S^2(2,3,4)$). In Section 6 we give another example of a
hyperbolic manifold admitting a reducible Dehn filling (i.e. the
connected sum of two lens spaces) and a Dehn filling of
icosahedral type (i.e. $S^2(2,3,5)$).

Throughout the paper, $S(a_1,a_2,\ldots,a_n)$ denotes a Seifert
fibered space over a surface $S$ with $n$ exceptional fibers of
orders $a_1,a_2,\ldots,a_n$, and $C(s,t)$ denotes the cable space
as defined in \cite[Section 3]{GL}.

\section{Cyclic Dehn fillings}
In this section, we show that there are hyperbolic manifolds which
admit a reducible Dehn filling and a finite cyclic Dehn filling at
distance 1.

The first example was given by Boyer and Zhang in \cite[Example
7.8]{BZ2}. In their example $M=W(6)$, which is obtained by Dehn
filling on one boundary component of the Whitehead link with slope
6, admits a reducible Dehn filling $M(1)$(=$L(3,1)\#L(2,1)$) and a
cyclic Dehn filling $M(1/0)$(=$L(6,1)$).

The first infinite family of hyperbolic manifolds admitting a
cyclic Dehn filling and a reducible Dehn filling was given by
Eudave-Mu\~{n}oz and Wu in \cite[Section 4]{EW}. In their
examples, the hyperbolic manifolds $M_{p}$, $p\geq2$ allow a
$0$-Dehn filling $M_{p}(0)$, which is the lens space
$L((p-1)(p+3)+1, p+3)$, and a $1/3$-Dehn filling $M_{p}(1/3)$,
which is the reducible manifold $L(3,1) \# L(2,1)$. We can observe
that Boyer and Zhang's manifold $W(6)$ belongs to this family
$M_p$ i.e. $W(6)=M_2$. Note that the reducible Dehn filling is
always $L(3,1) \# L(2,1)$ regardless of $p$. So it is natural to
try to find some examples realizing the connected sum of two lens
spaces of arbitrary orders.

\begin{figure}[t]
\includegraphics{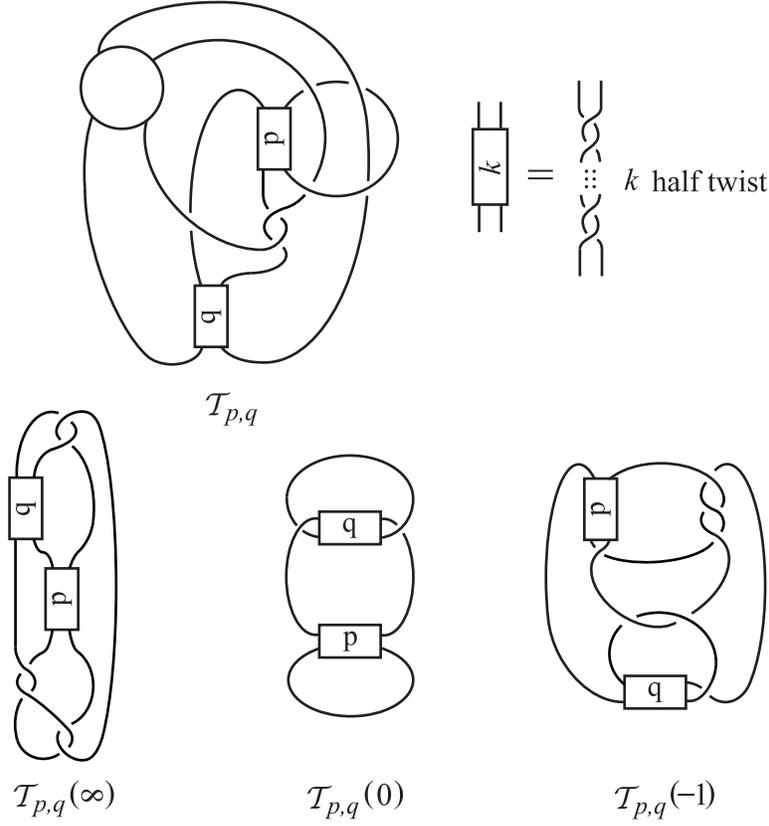}\caption{Cyclic Dehn fillings.}\label{fcyclic}
\end{figure}

We will construct hyperbolic manifolds which admit a Dehn filling
that is a connected sum of lens spaces of arbitrary orders and a
lens space Dehn filling. Consider the tangles $\T_{p,q}$
illustrated in Figure \ref{fcyclic} with $p\geq2$ and $q\geq4$.
Let $\T_{p,q}(r)$ be the link obtained by filling an $r$-rational
tangle to $\T_{p,q}$, $M_{p,q}$ the double branched cover of the
tangle $\T_{p,q}$ and $M_{p,q}(r)$ the double branched cover of
$S^3$ branched along $\T_{p,q}(r)$.

\begin{lemma}\label{cycliclemma}
The manifolds $M_{p,q}$ admit the following Dehn fillings.
\begin{itemize}
\item[(1)] $M_{p,q}(0)=L(p,1)\# L(q-2,1)$;

\item[(2)] $M_{p,q}(\infty)$ is the lens space
$L((3p+2)(-2q+1)+6,(3p+2)q-3)$;

\item[(3)] $M_{p,q}(-1)$ is an irreducible, toroidal and
non-Seifert fibered manifold.

\end{itemize}

\end{lemma}

\begin{proof}
The tangles $\T_{p,q}(0)$, $\T_{p,q}(\infty)$, $\T_{p,q}(-1)$ are
shown in Figure \ref{fcyclic}. The result now follows by taking
the double branched cover of $S^3$ branched along the
corresponding links.
\end{proof}

\begin{theorem}\label{cyclictheorem}
The manifolds $M_{p,q}$ are hyperbolic manifolds admitting two
Dehn fillings $M_{p,q}(r_1)$ and $M_{p,q}(r_2)$ such that
$M_{p,q}(r_1)$ is the connected sum of two lens spaces of
arbitrary orders, $M_{p,q}(r_2)$ is a lens space and
$\Delta(r_1,r_2)=1$.
\end{theorem}

\begin{proof}
Let $r_1=0$ and $r_2=\infty$. Then $\Delta(r_1,r_2)=1$ and by
Lemma \ref{cycliclemma}, $M(r_1)$ is the connected sum of two lens
spaces of arbitrary orders and $M(r_2)$ is a lens space. To
complete the proof, we need to show that $M_{p,q}$ is hyperbolic
i.e. $M_{p,q}$ is irreducible, $\partial$-irreducible, non-Seifert
fibered and atoroidal. However this follows from Lemmas
\ref{cycliclemma1}$-$\ref{cycliclemma2} below.

\end{proof}

\begin{lemma}\label{cycliclemma1}
$M_{p,q}$ is irreducible, $\partial$-irreducible and not Seifert
fibered.
\end{lemma}

\begin{proof}
First, we show that $M_{p,q}$ is irreducible. Suppose $M_{p,q}$ is
reducible. Then it contains a separating essential sphere since
$M_{p,q}(\infty)$ is a lens space. Therefore $M_{p,q}=X\#Y$ where
$Y$ contains the torus boundary of $M_{p,q}$ and $X\neq S^3$. Note
that if there are two slopes $r$ and $r'$ such that $M_{p,q}(r)$
and $M_{p,q}(r')$ are prime, then $M_{p,q}(r)=M_{p,q}(r')(=X)$
since $Y(r)$ and $Y(r')$ must be $S^3$. However, by Lemma
\ref{cycliclemma} $M_{p,q}(\infty)$ and $M_{p,q}(-1)$ are both
prime and distinct, which implies that $M_{p,q}$ is irreducible.

Suppose $M_{p,q}$ is a Seifert fibered space. Then $M_{p,q}(r)$ is
Seifert fibered for all but at most one $r$, for which
$M_{p,q}(r)$ is reducible. However $M_{p,q}(-1)$ is irreducible
and not Seifert fibered, which is a contradiction. Therefore
$M_{p,q}$ is not a Seifert fibered space.

Suppose $M_{p,q}$ is $\partial$-reducible. After
$\partial$-compression, the torus boundary becomes a sphere which
must bound a 3-ball since $M_{p,q}$ is irreducible. This implies
that $M_{p,q}$ is a solid torus, which is a Seifert fibered space,
a contradiction.
\end{proof}

To prove that $M_{p,q}$ is atoroidal, we need the following
lemmas.

\begin{lemma}\label{EWlemma}
Let $M(\neq T^2\times I)$ be an irreducible and
$\partial$-irreducible 3-manifold with a torus boundary component
$T_0$, and let $\alpha, \beta$ be slopes on $T_0$ with
$\Delta(\alpha,\beta)\geq 2$. Let $T_1$ be a torus in $\partial
M-T_0$ which is incompressible in $M$. If $T_1$ is compressible in
$M(\alpha)$ and $M(\beta)$, then $M$ is a cable space with cabling
slope $\gamma$ satisfying $\Delta(\alpha,
\gamma)=\Delta(\beta,\gamma)=1$.
\end{lemma}
\begin{proof}
This follows immediately from \cite[Theorems 2.4.3 and
2.4.4]{CGLS}.
\end{proof}

\begin{lemma}\label{generallemma}
Let $M(\neq T^2\times I)$ be an irreducible and
$\partial$-irreducible 3-manifold with a torus boundary component,
and let $\alpha, \beta$ be slopes on $T_0$ with
$\Delta(\alpha,\beta)=1$. Let $T_1$ be a torus in $\partial M-T_0$
which is incompressible in $M$. If $T_1$ is compressible in
$M(\alpha)$ and $M(\beta)$, then either
\begin{itemize}
\item[(1)]$M$ is a cable space $C(s,t)$, $t\geq2$ with cabling
slope $\alpha$, say, and thus $M(\alpha)=S^1 \times D^2 \# L(t,s)$
and $M(\beta)=S^1 \times D^2$. Furthermore if $r_\alpha, r_\beta$
are the slopes on $T_1$ corresponding to the meridians of the
solid tori of $M(\alpha), M(\beta)$ respectively, then
$\Delta(r_\alpha, r_\beta)\geq 2$; or

\item[(2)]$M(\alpha)$ and $M(\beta)$ are $S^1\times D^2$ and if
$r_\alpha, r_\beta$ are the slopes on $T_1$ corresponding to the
meridians of $M(\alpha), M(\beta)$ respectively, then
$\Delta(r_\alpha, r_\beta)\geq 4$.

\end{itemize}
\end{lemma}
\begin{proof}
It follows from \cite[Lemma 3.4]{EW} that either $M(\alpha)$ or
$M(\beta)$ must be irreducible. We assume without loss of
generality that $M(\beta)$ is irreducible. Since $T_1$ is
compressible in $M(\beta)$ i.e. $M(\beta)$ is
$\partial$-reducible, $M(\beta)=S^1 \times D^2$. Therefore
$\partial M$ consists of $T_0$ and $T_1$. Let $K_\beta$ be the
core of the Dehn filling solid torus. Then $M(\alpha)$ is obtained
by Dehn surgery on $K_\beta$ in the solid torus $M(\beta)$.

$(1)$ Assume that $M(\alpha)$ is reducible. It follows from
\cite[Theorem 1.1]{Ga} and \cite[Theorem 6.1]{Sch} that $M$ is a
cable space $C(s,t)$ with cabling slope $\alpha$. Hence
$M(\alpha)=S^1 \times D^2 \# L(t,s)$. Let $r_\alpha, r_\beta$ be
the slopes on $T_1$ corresponding to the meridians of the solid
tori summands of $M(\alpha), M(\beta)$ respectively. Since $M$ is
a cable space, we can apply \cite[Lemma 3.1]{GL} to get
$\Delta(r_\alpha, r_\beta)=|t|\Delta(\alpha, \beta)$. Since $M\neq
T^2 \times I$, $|t|\geq 2$. Therefore $\Delta(r_\alpha,
r_\beta)\geq 2$.

$(2)$ Assume that $M(\alpha)$ is irreducible. Since $T_1$ is
compressible in $M(\alpha)$ i.e. $M(\alpha)$ is
$\partial$-reducible, $M(\alpha)=S^1 \times D^2$. Then it follows
from \cite[Theorem 1.1]{Ga} that $M$ is the exterior of a braid in
a solid torus. Let $r_\alpha, r_\beta$ be the slopes on $T_1$
corresponding to the meridians of the solid tori $M(\alpha),
M(\beta)$ respectively. We can apply \cite[Lemma 3.3]{G2} to get
$\Delta(r_\alpha, r_\beta)\geq w^2$ where $w$ is the winding
number of $K_\beta$ in the solid torus $M(\beta)$. Since $M\neq
T^2 \times I$, $|w|\geq 2$. We are done.

\end{proof}

\begin{lemma}\label{cycliclemma2}
$M_{p,q}$ is atoroidal.
\end{lemma}

\begin{proof}
Suppose on the contrary that $M_{p,q}$ is not atoroidal i.e. it
contains an essential torus $F$. Note that $F$ is separating since
$M_{p,q}(\infty)$ is a lens space. Let $M_{p,q}=A\cup_{F} B$ where
$B$ contains $\partial M_{p,q}$. Since $M_{p,q}$ is irreducible
and $\partial$-irreducible, $A$ and $B$ are also irreducible and
$\partial$-irreducible. Lemma \ref{cycliclemma} implies that
$M_{p,q}(0)$ and $M_{p,q}(\infty)$ are atoroidal. Hence $F$ must
be compressible in both $B(0)$ and $B(\infty)$. Apply Lemma
\ref{generallemma} to $B$. Then there are two cases to consider.

\textit{Case 1: $B$ is a cable space $C(s,t)$}. First assume that
$\infty$ is the cabling slope. Then $B(\infty)=S^1\times D^2 \#
L(t,s)$ and $B(0)=S^1\times D^2$. Let $r_\infty, r_0$ be the
slopes on $F$ corresponding to the meridians of the solid tori of
$B(\infty), B(0)$ respectively. Then $\Delta(r_\infty, r_0)\geq
2$. Observe that $M(\infty)\cong A(r_\infty)\# L(t,s)$ and
$M(0)\cong A(r_0)$. Since $M_{p,q}(\infty)$ is a lens space and
$M_{p,q}(0)$ is the connected sum of two lens spaces,
$A(r_\infty)\cong S^3$ and $A(r_0)$ is reducible. Then
$\Delta(r_\infty, r_0)\leq 1$ by \cite{GLu3}, a contradiction.

Assume that $0$ is the cabling slope of $B=C(s,t)$. Then
$B(0)=S^1\times D^2 \# L(t,s)$ and $B(\infty)=S^1\times D^2$ and
$\Delta(r_0, r_\infty)\geq 2$. Moreover, $M(0)\cong A(r_0)\#
L(t,s)$ and $M(\infty)\cong A(r_\infty)$. By Lemma
\ref{cycliclemma}, $A(r_0)$ and $A(r_\infty)$ are lens spaces. By
the Cyclic Surgery Theorem of \cite{CGLS}, $A$ must be a Seifert
fibered space. Consider the $-1$-Dehn filling $M(-1)$, which is an
irreducible, toroidal and non-Seifert fibered 3-manifold. Since
$B=C(s,t)$ with the cabling slope $0$ and $\Delta(0,-1)=1$,
$B(-1)$ is a solid torus. Thus $M(-1) \cong A(r_{-1})$ where
$r_{-1}$ is the slope on $F$ corresponding to the meridian of the
solid torus $B(-1)$. Since $A$ is a Seifert fibered space,
$A(r_{-1})(\cong M(-1))$ is either a Seifert fibered space or a
reducible manifold, a contradiction.

\textit{Case 2: $B(\infty)$ and $B(0)$ are $S^1\times D^2$}. Then
$M(\infty)\cong A(r_\infty)$ and $M(0)\cong A(r_0)$. By Lemma
\ref{cycliclemma}, $A(r_\infty)$ is a lens space and $A(r_0)$ is
reducible. Also by Lemma \ref{generallemma} $\Delta(r_\infty,
r_0)\geq 4$. This is a contradiction to \cite[Theorem 1.2]{BZ2}.

\end{proof}

\section{Dihedral Dehn fillings}

In this section, we show that there are hyperbolic manifolds which
admit a reducible Dehn filling and a dihedral Dehn filling i.e. a
Dehn filling of type $S^2(2,2,n)$ at distance 1.

The first example was Boyer and Zhang's manifold $M=W(6)$ as
described in Section 2. The manifold $M=W(6)$ admits a Dehn
filling $M(2)$ which is $S^2(2,2,4)$, i.e. of \textit{dihedral}
type.

\begin{figure}[t]
\includegraphics{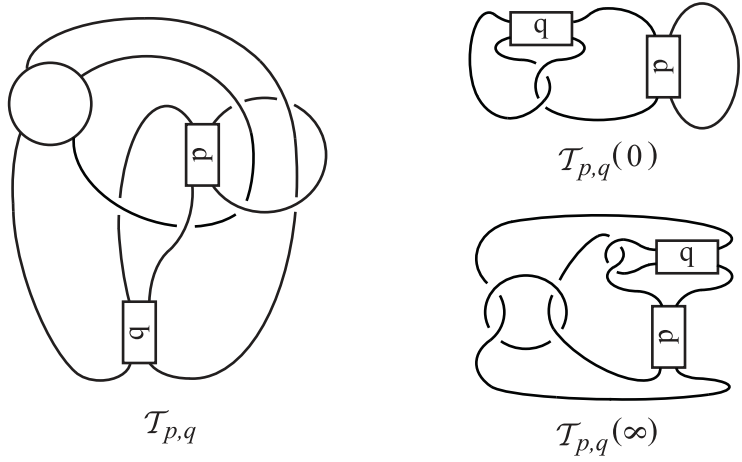}\caption{Dihedral Dehn fillings.}\label{Fdihedral}
\end{figure}

We present infinitely many examples of such hyperbolic manifolds.
Consider the tangles $\T_{p,q}$ with $p,q\geq 3$ illustrated in
Figure \ref{Fdihedral}. Then $0$-rational tangle filling on
$\T_{p,q}$, $\T_{p,q}(0)$, is a connected sum of two 2-bridge
links. Thus the double branched cover $M_{p,q}(0)$ is a connected
sum of two lens spaces. $\infty$-rational tangle filling on
$\T_{p,q}$, $\T_{p,q}(\infty)$, is a Montesinos link, whose double
branched cover is a Seifert fibered space over $S^2$ with three
exceptional fibers. More precisely, we get the following.

\begin{lemma}\label{Dlemma}
The manifolds $M_{p,q}$ admit the following Dehn fillings.
\begin{itemize}
\item[(1)] $M_{p,q}(0)=L(p,1)\# L(2q+1,1)$;

\item[(2)] $M_{p,q}(\infty)=S^2(2,2,n)$ where $n=2pq-p-2$.
\end{itemize}

\end{lemma}

The main theorem of this section is the following.
\begin{theorem}\label{Dtheorem}
The manifolds $M_{p,q}$ are hyperbolic manifolds, admitting two
Dehn fillings $M_{p,q}(r_1)$ and $M_{p,q}(r_2)$, such that
$M_{p,q}(r_1)$ is the connected sum of two lens spaces,
$M_{p,q}(r_2)$ is $S^2(2,2,n)$ and $\Delta(r_1,r_2)=1$.
\end{theorem}

\begin{proof}
This follows immediately from Lemmas \ref{Dlemma}, \ref{Dlemma1}
and \ref{Dlemma2}.
\end{proof}

\begin{lemma}\label{Dlemma1}
$M_{p,q}$ is irreducible, $\partial$-irreducible and not Seifert
fibered.
\end{lemma}

\begin{proof}
First, we show that $M_{p,q}$ is irreducible. Suppose $M_{p,q}$ is
reducible i.e. it contains an essential sphere $S$. Then $S$ must
be separating since $M_{p,q}(\infty)$ is irreducible. The sphere
$S$ decomposes $M_{p,q}$ as $X \# Y$ where $Y$ contains $\partial
M_{p,q}$. Then $M_{p,q}(\infty)=X\# Y(\infty)$. Since
$M_{p,q}(\infty)(=S^2(2,2,n))$ is irreducible, $Y(\infty)$ must be
$S^3$. Hence $X\cong S^2(2,2,n)$.

On the other hand, $0$-Dehn filling $M_{p,q}(0)$ is $X\# Y(0)$. If
$S$ is inessential in $M_{p,q}(0)$ i.e. $Y(0)$ is $S^3$, then
$M_{p,q}(0)\cong X$, which is the connected sum of two lens
spaces, a contradiction to $X=S^2(2,2,n)$. If $S$ is essential in
$M_{p,q}(0)$, then by the uniqueness of the prime decomposition of
a 3-manifold $X$ must be a lens space, a contradiction.

Secondly, we show that $M_{p,q}$ is not a Seifert fibered space.
Suppose $M_{p,q}$ is a Seifert fibered space. Since
$M_{p,q}(\infty)$ is $S^2(2,2,n)$, $M_{p,q}$ is either $D^2(a,b)$
, $D^2(a,b,c)$ or $M^2(c)$ where one of $a,b$ is 2 (let $a=2$) and
$M^2$ is a M\"obius band. Observe that $M_{p,q}(r)$ is Seifert
fibered for all but at most one $r$, for which $M_{p,q}(r)$ is
reducible. When $M_{p,q}(r)$ is reducible, $r$ corresponds to the
slope of the Seifert fiber of $M_{p,q}$. Lemma \ref{Dlemma} shows
that $0$ is the slope of the Seifert fiber of $M_{p,q}$.
Considering the fundamental group of $D^2(a,b)$, $D^2(a,b,c)$, or
$M^2(c)$ and Dehn filling with slope of the Seifert fiber, it
follows that $\pi_1(M_{p,q}(0))$ is isomorphic to either $\zed_a *
\zed_b$, $\zed_a * \zed_b * \zed_c$ or $\zed * \zed_c$ with $a=2$.
This is a contradiction since $\pi_1(M_{p,q}(0))=\zed_p *
\zed_{2q+1}$ where $p$ and $2q+1$ are greater than $2$.

$\partial$-irreducibility of $M_{p,q}$ follows from irreducibility
and the fact that $M_{p,q}$ is not Seifert fibered.

\end{proof}

\begin{lemma}\label{Dlemma2}
$M_{p,q}$ is atoroidal.
\end{lemma}

\begin{proof}
Suppose on the contrary that $M_{p,q}$ is toroidal. We consider
the tangle $\T_p$ with two tangle spheres $\S_0, \S_1$ depicted in
Figure \ref{Fdihedral1}. Let $N_p, \partial_0 N_p, \partial_1 N_p$
be the double branched covers of $\T_p, \S_0, \S_1$ respectively.
Then observe that the tangle $\T_{p,q}$ in Figure \ref{Fdihedral}
is obtained from $\T_p$ by $1/q$-rational tangle filling on
$\S_1$. From the viewpoint of the double branched cover, $M_{p,q}$
is obtained from $N_p$ by $1/q$-Dehn filling on $\partial_1 N_p$
i.e. $M_{p,q}=N_p(1/q)$. Hence $N_p(1/q)$ is toroidal.

\begin{figure}[t]
\includegraphics{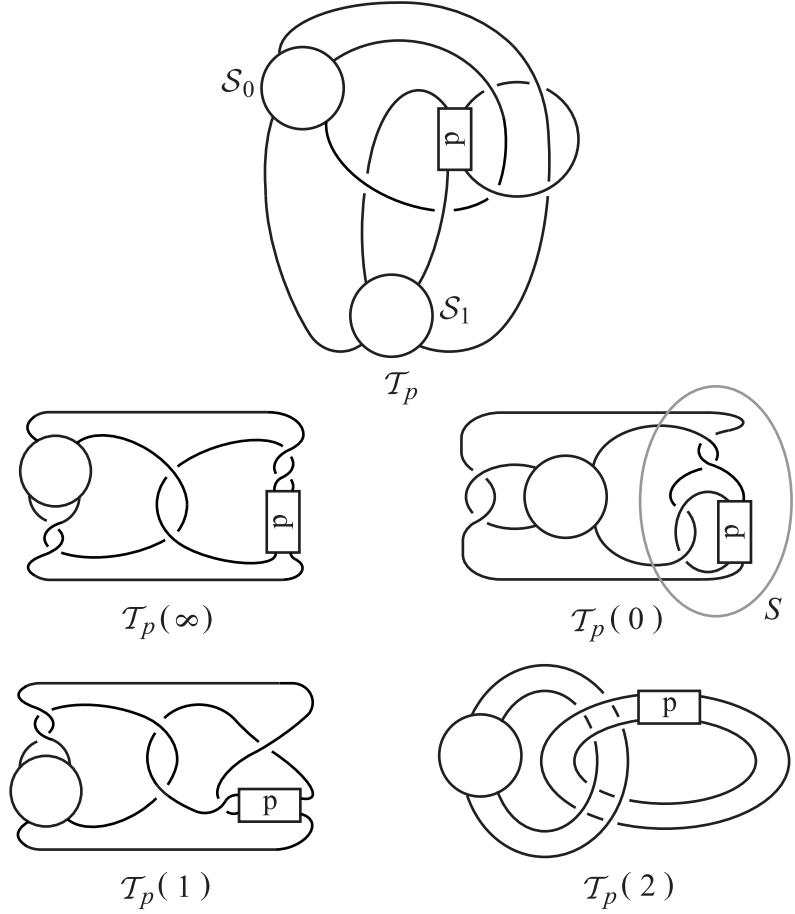}\caption{Rational tangle fillings on $\T_p$
obtained from $\T_{p,q}$ by removing $1/q$-rational
tangle.}\label{Fdihedral1}
\end{figure}

We perform rational tangle fillings on $\T_p$ along $\S_1$.
Several rational tangle fillings are shown in Figure
\ref{Fdihedral1}. Then we get the corresponding double branched
covers as follows; $N_p(\infty)\cong D^2(2, p+2)$, $N_p(0)\cong
C(1,2)\cup_R D^2(2,p),$ $N_p(1)\cong D^2(2, p-2)$ and $N_p(2)\cong
Z\times S^1$ where $R$ is a torus corresponding to the Conway
sphere $S$ as shown in Figure \ref{Fdihedral1} and $Z$ is a
once-punctured torus. Note that $N_p(0)$ is not Seifert fibered
and $N_p(2)$ contains a non-separating essential torus.

\begin{claim}
$N_p$ is hyperbolic.
\end{claim}
\begin{proof}
Since $N_p(\infty)$ and $N_p(1)$ are distinct prime manifolds,
$N_p$ is irreducible. Also, $N_p$ is not Seifert fibered because
$N_p(1/q)$ is an irreducible non-Seifert fibered space by Lemma
\ref{Dlemma1}. By irreducibility and the fact that it is not
Seifert fibered, $N_p$ is $\partial$-irreducible. To complete the
proof, we need only to show that $N_p$ is atoroidal.

Suppose $N_p$ contains an essential torus $F$. Then $F$ must be
separating since $N_p(\infty)$ does not contain a non-separating
torus or sphere. Let $N_p=A \cup_F B$ with the filling torus
$\partial_1 N_p \subseteq B$.

Recall that $N_p(\infty)\cong D^2(2, p+2), N_p(1)\cong D^2(2,
p-2)$, $N_p(0)\cong C(1,2)\cup_R D^2(2,p)$ and $N_p(2)\cong
Z\times S^1$ contains a non-separating essential torus. Thus $F$
must be compressible in $B(\infty), B(1)$ and $B(2)$. Apply Lemma
\ref{generallemma} to $B$ with slopes $\infty, 1,2$. Since
$N_p(\infty), N_p(1)$ and $N(2)$ don't have a lens space summand,
$\infty, 1$ and $2$ can't be cabling slope, which implies that
$B(\infty), B(1)$ and $B(2)$ are solid tori. Therefore $F$
separates the two boundaries $\partial_0 N_p$, $\partial_1 N_p$ of
$N_p$. In other words, $A$ has the boundary $\partial_0 N_p$ and
$B$ has the boundary $\partial_1 N_p$.

Suppose $F$ is incompressible in $N_p(0)\cong C(1,2)\cup_R
D^2(2,p)$. Then $F$ is isotopic to $R$. Thus $A$ is homeomorphic
to either $C(1,2)$ or $D^2(2,p)$. However $A$ has the boundary
$\partial_0 N_p$ and therefore $A\cong C(1,2)$. Since $B(2)$ is a
solid torus, $N_p(2)$ can be obtained from $C(1,2)$ by Dehn
filling on $F$. Hence $N_p(2)$ does not contain a non-separating
torus, a contradiction. It follows that $F$ is compressible in
$N_p(0)$ and thus in $B(0)$.

We have shown that $F$ is compressible in $B(\infty), B(0), B(1)$
and $B(2)$. Since $\Delta(0,2)=2$, it follows from Lemma
\ref{EWlemma} that $B$ is a cable space with cabling slope either
$\infty$ or $1$. This implies that either $N_p(\infty)$ or
$N_p(1)$ has a lens space summand, which is a contradiction. This
completes the proof of the claim.
\end{proof}

The claim says that $N_p$ is hyperbolic. However $N_p$ admits the
two toroidal Dehn fillings $N_p(2), N_p(1/q)$, and $\Delta(2,
1/q)=|2q-1|\geq 5$ since $q\geq 3$. Since $N_p$ has two boundary
components, \cite[Theorem 1.1]{GW} implies that $N_p$ is
homeomorphic to the exterior of the Whitehead sister link and in
the two toroidal Dehn fillings every essential torus is
separating. This is a contradiction because $N_p(2)$ contains a
non-separating essential torus.
\end{proof}
%%%%%%%%%%%%%%%%%%%%%%%%%%%%%%%%%%%%%%%%%%%%%%%%%%%%%%%%%%%%%%%%%%%%%%%%%%%%%%

\section{Tetrahedral Dehn fillings}

\begin{figure}[t]
\includegraphics{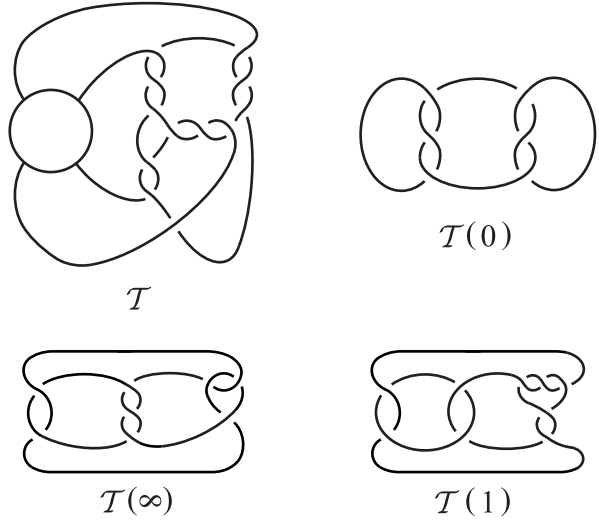}\caption{Tetrahedral Dehn fillings.}\label{Ft}
\end{figure}

In this section, we describe a hyperbolic manifold which admits a
reducible Dehn filling and a tetrahedral Dehn filling, i.e. an
$S^2(2,3,3)$-Dehn filling.

Let $\T$ be the tangle depicted in Figure \ref{Ft}. We perform 0,
1 and $\infty$-rational tangle fillings on $\T$. See Figure
\ref{Ft}. Then $\T(0)$ is the connected sum of two trefoil knots
and $\T(1), \T(\infty)$ are Montesinos links. Let $M$ be the
double branched cover of the tangle $\T$. Considering double
branched covers, we have the following lemma.

\begin{lemma}\label{Tlemma1}
The manifold $M$ admits the following Dehn fillings.
\begin{itemize}
\item[(1)] $M(0)=L(3,1)\#L(3,1)$;

\item[(2)] $M(\infty)=S^2(2, 3, 3)$;

\item[(3)] $M(1)=S^2(2, 2, 7)$.
\end{itemize}

\end{lemma}

\begin{lemma}\label{Tlemma2}
$M$ is irreducible, $\partial$-irreducible and not Seifert
fibered.
\end{lemma}

\begin{proof}
Since $M(\infty)$ and $M(1)$ are distinct prime manifolds, $M$ is
irreducible.

Next, we show that $M$ is not Seifert fibered. Suppose on the
contrary that $M$ is a Seifert fibered space. Since $\infty$-Dehn
filling $M(\infty)$ is $S^2(2, 3, 3)$, $M$ is $D^2(a,b)$ or
$D^2(a,b,c)$ where $a$, say, must be 3. Then $M(1)=S^2(2,2,7)$ is
obtained by Dehn filling on $D^2(a,b)$ or $D^2(a,b,c)$. This is
impossible since $a$ is 3.

The $\partial$-irreducibility follows from the above two facts
about $M$.

\end{proof}

\begin{theorem} \label{Ttheorem}
The manifold $M$ is a hyperbolic manifold admitting two Dehn
fillings $M(r_1)$ and $M(r_2)$ such that $M(r_1)$ is
$L(3,1)\#L(3,1)$, $M(r_2)$ is $S^2(2,3,3)$ and $\Delta(r_1,
r_2)=1$.
\end{theorem}

\begin{proof}
By Lemmas \ref{Tlemma1}, \ref{Tlemma2}, we need only to show that
$M$ is atoroidal.

Suppose $M$ is toroidal. Since $M(\infty)$ does not contain a
non-separating torus or sphere, any essential torus in $M$ must be
separating. Let $F$ be an innermost essential torus in $M$, i.e.
$M=A\cup_F B$ where $\partial M \subseteq B$ and $A$ is atoroidal.
Then $A$ and $B$ are irreducible and $\partial$-irreducible since
$M$ is. By Lemma \ref{Tlemma1}, $F$ must be compressible in
$M(\infty)$ and $M(1)$. Thus $F$ is compressible in $B(\infty)$
and $B(1)$. Apply Lemma \ref{generallemma} to $B$ with slopes
$\infty$ and $1$. Since $M(\infty)$ and $M(1)$ don't have a lens
space summand, the first case of Lemma \ref{generallemma} can't
occur. Thus $B(\infty)$ and $B(1)$ are solid tori and if we let
$r_\infty$ and $r_1$ be the slopes on $F$ corresponding to the
meridians of $B(\infty)$ and $B(1)$ respectively, then
$\Delta(r_\infty, r_1)\geq 4$. Also $M(\infty) \cong A(r_\infty)$
and $M(1) \cong A(r_1)$. In other words, $A(r_\infty)=S^2(2, 3,
3)$ and $A(r_1)=S^2(2, 2, 7)$.

Recall that $A$ is irreducible, $\partial$-irreducible and
atoroidal. Furthermore, it is easy to see by applying an argument
similar to that of the second paragraph of Lemma \ref{Tlemma2}
that $A$ is not Seifert fibered. Therefore $A$ is hyperbolic.
Since $A(r_\infty)=S^2(2, 3, 3)$ and $A(r_1)=S^2(2, 2, 7)$ with
$\Delta(r_\infty, r_1)\geq 4$, $A$ admits two finite Dehn fillings
with distance $\geq 4$. This is a contradiction to \cite[Theorem
1.1]{BZ3}.

\end{proof}

\textbf{Remark.} Lemma \ref{Tlemma1} shows that the manifold $M$
in Theorem \ref{Ttheorem} is also an example of a hyperbolic
manifold admitting a \textit{dihedral} finite Dehn filling.

%%%%%%%%%%%%%%%%%%%%%%%%%%%%%%%%%%%%%%%%%%%%%%%%%%%%%%%%%%%%%%%%%%%%%%%%%%%%%%%%%%%%%%%%%%%

\section{Octahedral Dehn fillings}
In this section, we present hyperbolic manifolds admitting a
reducible Dehn filling and an octahedral Dehn filling i.e. an
$S^2(2,3,4)$-Dehn filling.

\begin{figure}[t]
\includegraphics{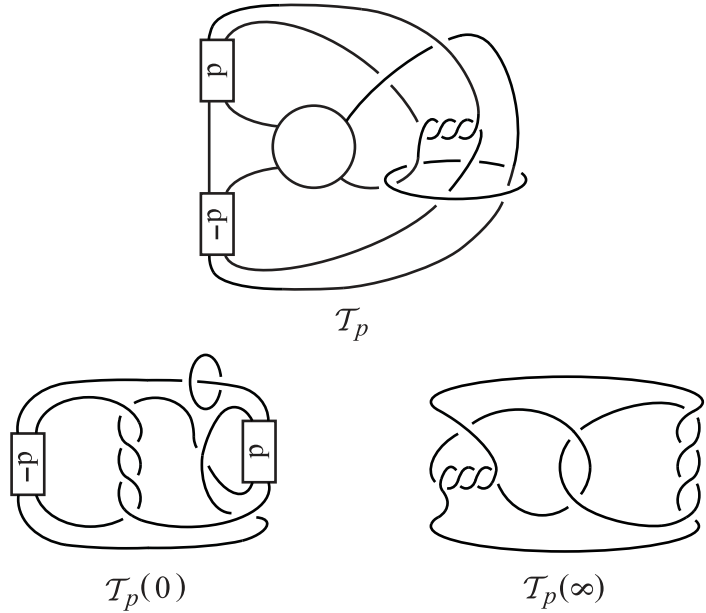}\caption{Octahedral Dehn fillings.}\label{Fo}
\end{figure}

We consider the tangle $\T_p$, $p\geq 3$ illustrated in Figure
\ref{Fo}. Then $0$-rational tangle filling $\T_p(0)$ gives the
connected sum of the Hopf link and a Montesinos link.
$\infty$-rational tangle filling $\T_p(\infty)$ gives rise to a
Montesinos link. Let $M_p$ be the double branched cover of the
tangle $\T_p$. Then the following lemma follows immediately by
considering the double branched covers of the above tangles or
links.

\begin{lemma}\label{Olemma}
The manifolds $M_p$ admit the following Dehn fillings.
\begin{itemize}
\item[(1)] $M_p(0)=L(2,1)\# S^2(4, p, 2p+1)$;

\item[(2)] $M_p(\infty)=S^2(2,3,4)$.
\end{itemize}

\end{lemma}

\begin{theorem}\label{Otheorem}
The manifolds $M_p$ are hyperbolic manifolds admitting two Dehn
fillings $M_p(r_1)$ and $M_p(r_2)$ such that $M_p(r_1)$ is $L(2,1)
\# S^2(4,$ $p,2p+1)$, $M_p(r_2)$ is $S^2(2,3,4)$, and
$\Delta(r_1,r_2)=1$.
\end{theorem}

\begin{proof}
This follows immediately from Lemmas \ref{Olemma}, \ref{Olemma1}
and \ref{Olemma2}.
\end{proof}

\begin{lemma}\label{Olemma1}
$M_p$ is irreducible, $\partial$-irreducible and not Seifert
fibered.
\end{lemma}

\begin{proof}
First, we show that $M_p$ is irreducible.

Suppose $M_p$ is reducible. Then it contains an essential
separating sphere $S$, inducing a decomposition $M_p=X\# Y$ where
$Y$ contains $\partial M_p$. $\infty$-Dehn filling $M_p(\infty)$
is $X\# Y(\infty)$. By the irreducibility of
$M_p(\infty)(=S^2(2,3,4))$, $Y(\infty)$ must be $S^3$ and thus
$X=S^2(2,3,4)$.

$0$-Dehn filling $M_p(0)$ is $X\# Y(0)$. $S$ must be essential in
$M_p(0)$. Otherwise $Y(0)$ is $S^3$ and $M_{p,q}(0)=X$ is
reducible, a contradiction to $X=S^2(2,3,4)$. By the uniqueness of
the prime decomposition of a 3-manifold and Lemma \ref{Olemma},
$X$ must be either $L(2,1)$ or $S^2(4,p,2p+1)$. This is a
contradiction because $X=S^2(2,3,4)$.

Suppose $M_p$ is a Seifert fibered space.  Since $M_p(\infty)$ is
$S^2(2,3,4)$, $M_p$ is either $D^2(a,b)$ or $D^2(a,b,c)$ where one
of $a,b,c$, say $a$, is either $2,3,$ or $4$. $M_p$ admits a
non-Seifert fibered Dehn filling $M_p(r)$ for at most one slope
$r$, for which $M_p(r)$ is reducible. Furthermore such an $r$ is
the slope of the Seifert fiber of $M_p$. Since $M_p(0)$ is
reducible, $0$ is the slope of the Seifert fiber of $M_p$. A
fundamental group argument shows that $\pi_1(M_p(0))$ is
isomorphic to $\mathbb Z_a
* \mathbb Z_b$ or $\mathbb Z_a * \mathbb Z_b * \mathbb Z_c$, where $a=2,3$ or $4$.
This is a contradiction since $\pi_1(M_p(0))=\mathbb Z_2 * G$
where $G=\pi_1(S^2(4,p,2p+1))$

$\partial$-irreducibility of $M_p$ follows from irreducibility and
the fact that $M_p$ is non-Seifert fibered.

\end{proof}

\begin{lemma}\label{Olemma2}
$M_p$ is atoroidal.
\end{lemma}

\begin{proof}
Assuming the contrary, let $F$ be an essential torus in $M_p$.
Since $M_p(\infty)$ does not contain a non-separating torus or
sphere, $F$ must be separating. Let $M_p=A\cup_{F} B$ where $B$
contains $\partial M_p$. We may choose $F$ to be innermost, so
that $A$ is atoroidal.

Since $M_p$ is irreducible and $\partial$-irreducible, $A$ and $B$
are also irreducible and $\partial$-irreducible. Observe by Lemma
\ref{Olemma} that $M_p(0)$ and $M_p(\infty)$ are atoroidal, which
implies that $F$ must be compressible in both $B(0)$ and
$B(\infty)$. Then by Lemma \ref{generallemma}, there are two cases
to consider.

\textit{Case 1: $B$ is a cable space $C(s,t)$, $t\geq 2$ and
either $0$ or $\infty$ is the cabling slope.} Suppose that
$\infty$ is the cabling slope of $B$. Then $B(\infty)=S^1\times
D^2 \# L(t,s)$ and thus $M_p(\infty)$ has a lens space summand,
which is a contradiction to Lemma \ref{Olemma}. Hence the cabling
slope must be $0$. Then $B(0)=S^1\times D^2 \# L(t,s)$ and
$M_p(0)$ has a lens space summand $L(t,s)$. Lemma \ref{Olemma}
implies that $L(t,s)\cong L(2,1)$ i.e. $t=2$.

\textit{Case 2: $B(\infty)$ and $B(0)$ are $S^1\times D^2$}. Let
$r_\infty, r_0$ be the slopes on $F$ corresponding to the
meridians of the solid tori $B(\infty), B(0)$ respectively. Then
$M_p(\infty)\cong A(r_\infty)$ and $M_p(0)\cong A(r_0)$.
Correspondingly, $A(r_\infty)=S^2(2,3,4)$ and $A(r_0)$ is
reducible. Also $\Delta(r_\infty, r_0)\geq 4$ by Lemma
\ref{generallemma}. Recall that $A$ is an irreducible,
$\partial$-irreducible and atoroidal 3-manifold with torus
boundary $F$. However, the non-Seifert fiberedness of $A$ follows
from an argument similar to that proving that $M_p$ is not a
Seifert fibered space in Lemma \ref{Olemma1}. This implies that
$A$ is hyperbolic, and admits a reducible Dehn filling and a
finite Dehn filling with distance greater than 1. This is a
contradiction to \cite[Theorem 1]{BGZ}.

We have shown from Cases $1,2$ that $B$ is the cable space
$C(s,2)$ with cabling slope $0$. Hence $M_p=A \cup_F C(s,2)$.

\begin{figure}[t]
\includegraphics{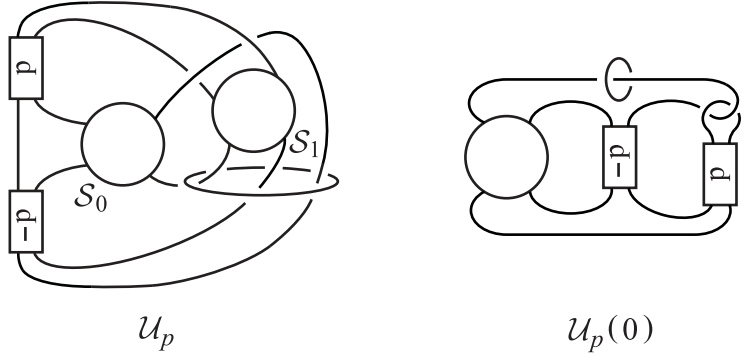}\caption{Rational tangle fillings on $\U_p$ obtained from $\T_p$ by removing 4-rational tangle.}\label{Fo1}
\end{figure}

Now we consider the tangle $\U_p$ with two tangle spheres $\S_0,
\S_1$ shown in Figure \ref{Fo1}. Let $N_p, \partial_0 N_p,
\partial_1 N_p$ be the double branched covers of $\U_p, \S_0, \S_1$
respectively. Then observe that the tangle $\T_p$ in Figure
\ref{Fo} is obtained from $\U_p$ by $4$-rational tangle filling on
$\S_1$. As double branched covers, $M_p$ is obtained from $N_p$ by
$4$-Dehn filling on $\partial_1 N_p$ i.e. $M_p=N_p(4)$. Hence
$N_p(4)=A \cup_F C(s,2)$, which is toroidal.

\begin{claim}\label{separate}
$F$ separates the two boundaries $\partial_0 N_p$ and $\partial_1
N_p$ of $N_p$. Therefore $A$ has the boundary $\partial_1 N_p$.
\end{claim}

\begin{proof}
Since $B$ is the cable space $C(s,2)$ with cabling slope $0$,
$B(0)=S^1\times D^2 \# L(2,s)$. Let $V$ be the solid torus summand
of $B(0)$ and $J$ the core of $V$. Note that $\partial V=F$ since
$\partial B(0)=F$. Consider $M_p(0)=L(2,1)\# S^2(4, p, 2p+1)$.
Since $M_p(0)=A\cup_F B(0)$ and $B(0)=V \# L(2,s)$, it follows
that $S^2(4, p, 2p+1)$ is obtained from $A$ by attaching $V$ along
$F$ i.e. $S^2(4, p, 2p+1)=A \cup_F V$. Therefore $S^2(4, p, 2p+1)$
contains $V$.

Recall that the tangle $\U_p$ is obtained from the tangle $\T_p$
by removing the 4-rational tangle. Let $R$ be the 4-rational
tangle. Let $W$ be the solid torus in $M_p$ which is the double
branched cover of $R$, and $K$ the core of $W$. Note that
$\partial R=\S_1$. Therefore $\partial W= \partial_1 N_p$. Recall
from Figure \ref{Fo} that $\T_p(0)$ is the connected sum of the
Hopf link and a Montesinos link. Then $S^2(4, p, 2p+1)$ summand in
$M_p(0)$ is the double branched cover of the Montesinos link in
$\T_p(0)$. Furthermore we can observe from Figure \ref{Fo} that
$-1/4$-rational tangle in the Montesinos link comes directly from
the 4-rational tangle $R$. Therefore $S^2(4, p, 2p+1)$ contains
$W$ since $W$ is the double branched cover of $R$.

We have seen that $S^2(4, p, 2p+1)$ contains $V$ and $W$ with
$\partial V=F$ and $\partial W=\partial_1 N_p$. Perturbing the
cores $J$ and $K$ of $V$ and $W$, we may assume that $K$ doesn't
intersect $J$ in $S^2(4, p, 2p+1)$. In other words, $W$ doesn't
intersect $V$. Thus $\partial_1 N_p$ doesn't intersect $V$. Since
$S^2(4, p, 2p+1)=A \cup_F V$, $\partial_1 N_p$ is contained in
$A$, as desired.

\end{proof}

We perform a $0$-rational tangle filling on $\U_p$ along $\S_1$,
which gives the connected sum of the Hopf link and some Montesinos
tangle as described in Figure \ref{Fo1}. Then the corresponding
double branched cover $N_p(0)$ is $L(2,1)\# D^2(p,$ $2p+1)$. We
claim that $N_p$ is hyperbolic. As long as the claim holds, the
hyperbolic manifold $N_p$ admits the reducible Dehn filling
$N_p(0)$ and the toroidal Dehn filling $N_p(4)$, and
$\Delta(0,4)=4$. This is a contradiction to \cite[Theorem 1.1]{O}
or \cite[Theorem 1]{W}.

\begin{claim}
$N_p$ is hyperbolic.
\end{claim}

\begin{proof}
$N_p$ admits two exceptional Dehn fillings $N_p(0)= L(2,1)\#
D^2(p$, $2p+1)$ and $N_p(4)=A \cup_F C(s,2)$, which is toroidal.
Note that by Lemma \ref{Olemma1} $N_p(4)$ is irreducible,
$\partial$-irreducible and non-Seifert fibered.

Suppose that $N_p$ is not irreducible. In other words, there is an
essential sphere $S$ in $N_p$. $S$ must be separating since
$N_p(4)$ does not contain a non-separating sphere, so it induces a
decomposition $N_p= X\# Y$ with $\partial_1 N_p \subseteq Y$. Then
$N_p(0)=X \# Y(0)$ and $N_p(4)= X \# Y(4)$ . Since $N_p(4)$ is
irreducible, $Y(4)$ is $S^3$ and thus $N_p(4)\cong X$. Hence $X$
is irreducible and toroidal. On the other hand, if $S$ is
inessential in $N_p(0)$, then $Y(0)$ must be $S^3$, which implies
that $X(\cong N_p(0))$ is reducible, a contradiction. If $S$ is
essential in $N_p(0)$, then $X$ must be either $L(2,1)$ or $D^2(p,
2p+1)$, neither of which is toroidal, a contradiction.

$N_p$ is not Seifert fibered since it has two non-Seifert fibered
Dehn fillings $N_p(0)$ and $N_p(4)$. By irreducibility and the
fact that it is non-Seifert fibered, $N_p$ is
$\partial$-irreducible. To complete the proof, we need only to
show that $N_p$ is atoroidal.

Suppose $N_p$ contains an essential torus $G$. Then $G$ must be
separating since $N_p(0)$ does not contain a non-separating torus
or sphere. Let $N_p=D \cup_G E$ where the filling torus
$\partial_1 N_p \subseteq E$ and $G$ is chosen to be innermost
with respect to $\partial_1 N_p$, so that $D$ is atoroidal. Then
$D$ and $E$ are irreducible and $\partial$-irreducible since $N_p$
is. $G$ is compressible in $N_p(0)$ and thus in $E(0)$ since
$N_p(0)$ is atoroidal. Hence $E(0)=S^1\times D^2 \# W$ for some
3-manifold $W$. Let $s_0$ be the slope on $G$ corresponding to the
meridian of the solid torus summand of $E(0)$. Then
$N_p(0)=D(s_0)\# W$.

Suppose $G$ is compressible in $N_p(4)$. By Lemma \ref{EWlemma}
$E$ is a cable space $C(u,v)$ with $v\geq 2$ with cabling slope
$\infty$. Since $\Delta(\infty, 0)=\Delta(\infty, 4)=1$, $E(0)$
and $E(4)$ are solid tori $S^1\times D^2$. Let $s_4$ be the slope
on $G$ corresponding to the meridian of $E(4)$. Then
$N_p(0)=D(s_0)$ and $N_p(4)=D(s_4)$, which are reducible and
toroidal respectively. Moreover by \cite[Lemma 3.1]{GL},
$\Delta(s_0, s_4)=v^2\Delta(0,4)$ and thus $\Delta(s_0, s_4)\geq
16$. Recall that $D$ is irreducible, $\partial$-irreducible and
atoroidal. However $D$ is also not Seifert fibered. Otherwise, $D$
admits only one non-Seifert fibered Dehn filling, which is a
contradiction since $D(s_0)$ and $D(s_4)$ are not Seifert fibered.
Hence we have obtained a hyperbolic manifold $D$ which admits a
reducible Dehn filling $D(s_0)$ and a toroidal Dehn filling
$D(s_4)$ with distance $\geq 16$. This is a contradiction to
\cite[Theorem 1.1]{O} or \cite[Theorem 1]{W}.

Suppose $G$ is incompressible in $N_p(4)$. Since $N_p(4)$ is not a
Seifert fibered space, there is only one essential torus in
$N_p(4)=A \cup_F C(s,2)$ up to isotopy. Thus $G$ is isotopic to
$F$. Since $N_p(4)=D\cup_G E(4)$, $E(4)$ is homeomorphic to either
$A$ or $C(s,2)$. However Claim \ref{separate} says that $A$
contains $\partial_1 N_p$. Since $E(4)$ also contains $\partial_1
N_p$, $E(4)$ is homeomorphic to $A$. Therefore $D$ is homeomorphic
to $C(s,2)$.

Recall that $N_p(0)=D(s_0)\# W$. Since $D=C(s,2)$, $D(s_0)$ is
either $D^2(2,a)$ or $S^1\times D^2\# L(2,s)$ for some positive
integer $a$. Thus $N_p(0)= D^2(2,a)\# W$ or $ S^1\times D^2\#
L(2,s)\# W$. However either case is impossible since $N_p(0)=
L(2,1)\# D^2(p$, $2p+1)$ where $p\geq 3$. This completes the proof
of the claim.
\end{proof}

The claim completes the proof of the lemma.
\end{proof}

%%%%%%%%%%%%%%%%%%%%%%%%%%%%%%%%%%%%%%%%%%%%%%%%%%%%%%%%%%%%%%%%%%%%%%%%%%%

\section{Icosahedral Dehn fillings}

\begin{figure}[t]
\includegraphics{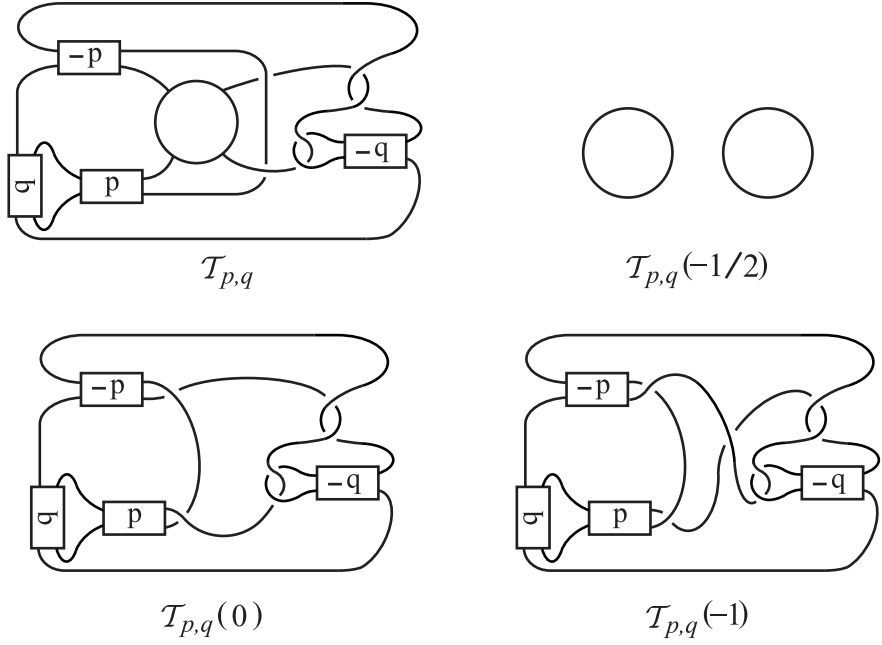}\caption{Icosahedral Dehn fillings.}\label{Fti}
\end{figure}

In this section, we present hyperbolic manifolds which admit a
reducible Dehn filling and a Dehn filling of type $S^2(2,3,5)$ at
distance 1.

Consider the tangles $\T_{p,q}$ shown in Figure \ref{Fti}, where
$p, q$ are integers satisfying $p\neq 0, \pm 1, q\neq 0, (p,q)\neq
(\pm 2, \pm 1)$. The tangles $\T_{p,q}$ were constructed by Lee
\cite[Section 4]{Lee2}. He used the tangles $\T_{p,q}$ to show
that there are hyperbolic 3-manifolds admitting an $S^1\times S^2$
Dehn filling and a toroidal Dehn filling with distance 2. (More
explicitly, the double branched covers of $\T_{p,q}(-1/2)$ and $
\T_{p,q}(\infty)$ are $S^1\times S^2$ and a toroidal 3-manifold
respectively.) Let $M_{p,q}$ be the double branched cover of
$\T_{p,q}$. Then $M_{p,q}$ with $p,q$ as above is hyperbolic by
\cite{Lee2}.

We consider several rational tangle fillings on $\T_{p,q}$,
$\T_{p,q}(-1/2)$, $\T_{p,q}(0)$, $\T_{p,q}(-1)$, which are
depicted in Figure \ref{Fti}. Then the following lemma follows
from Figure \ref{Fti}, considering double branched covers.

\begin{lemma}\label{Il}
The manifold $M_{p,q}$ admits the following Dehn fillings.
\begin{itemize}
\item[(1)] $M_{p,q}(-1/2)=S^1\times S^2$;

\item[(2)] $M_{p,q}(0)=S^2(|p-1|, |2q-1|, |pq+q-1|)$;

\item[(3)] $M_{p,q}(-1)=S^2(|p+1|, |2q+1|, |pq-q-1|)$.
\end{itemize}
\end{lemma}

\begin{theorem}\label{Itheorem}
The manifolds $M_{p,q}$, $(p,q)=(\pm3,\mp1)$, $(\pm4,\pm1)$, are
hyperbolic manifolds admitting two Dehn fillings $M_{p,q}(r_1)$
and $M_{p,q}(r_2)$, such that $M_{p,q}(r_1)$ is $S^1\times S^2$
i.e. reducible, $M_{p,q}(r_2)=S^2(2,3,5)$ and $\Delta(r_1,r_2)=1$.
\end{theorem}
\begin{proof}
As mentioned above, the manifolds $M_{p,q}$ are hyperbolic by
\cite{Lee2}.

If $(p,q)=(3,-1), (-4,-1)$, then $M_{p,q}(0)=S^2(2,3,5)$. If
$(p,q)=(-3,1), (4,1)$, then $M_{p,q}(-1)=S^2(2,3,5)$. Also
$\Delta(0,-1/2)=1$ and $\Delta(-1,$ $-1/2)=1$. We are done.
\end{proof}

\begin{figure}[t]
\includegraphics{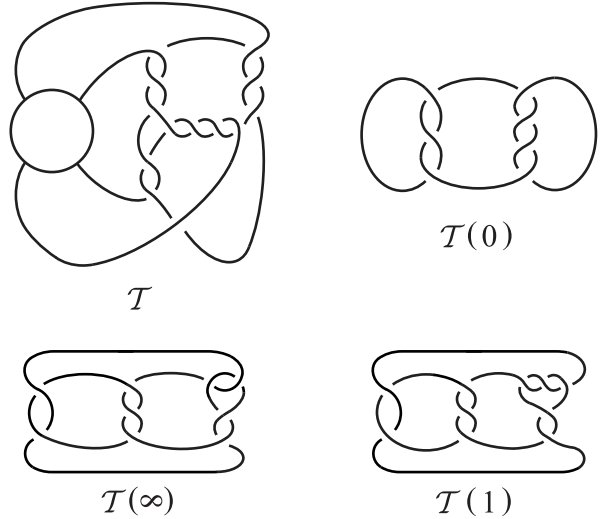}\caption{Another icosahedral Dehn fillings.}\label{Fti1}
\end{figure}

We present another hyperbolic manifold admitting an icosahedral
Dehn filling. Consider the tangle $\T$ and the three rational
tangle fillings $\T(0)$, $\T(\infty)$, $\T(1)$ shown in Figure
\ref{Fti1}. $\T(0)$ is the connected sum of the trefoil knot and
the $(2,4)$ torus link, while $\T(\infty)$ and $\T(1)$ are
Montesinos links. Let $M$ be the double branched cover of $\T$.
Considering the double branched covers of these, we have the
following lemma.

\begin{lemma} \label{I2}
The manifold $M$ admits three Dehn fillings as follows;
\begin{itemize}
\item[(1)] $M(0)=L(3,1)\# L(4,1)$;

\item[(2)] $M(\infty)=S^2(2,3,5)$;

\item[(3)] $M(1)=S^2(2,3,7)$.
\end{itemize}
\end{lemma}

\begin{theorem} \label{I3}
The manifold $M$ is a hyperbolic manifold. Thus $M$ is another
example of a hyperbolic manifold having a reducible and an
icosahedral Dehn filling.
\end{theorem}

\begin{proof}
Apply the same argument of Lemma \ref{Tlemma2} for irreducibility
and $\partial$-irreducibility of $M$, replacing $S^2(2,3,3)$,
$S^2(2,2,7)$ by $S^2(2,3,5)$, $S^2(2,$ $3,7)$ respectively.

Assume that $M$ is Seifert fibered. Then $M$ must be homeomorphic
to $D^2(2,3)$ by Lemma \ref{I2} $(2), (3)$. Observe that $M(r)$ is
Seifert fibered for all but at most one slope $r$, for which
$M(r)$ is reducible. When $M(r)$ is reducible, $r$ is the slope of
the Seifert fiber of $M$. Hence Lemma \ref{I2} implies that
$M(0)=L(3,1)\# L(4,1)$ is obtained by Dehn filling on $M=D^2(2,3)$
with the slope 0 corresponding to the Seifert fiber of $M$. This
is impossible by a fundamental group argument.

In order to prove that $M$ is hyperbolic, it remains to show that
$M$ is atoroidal. Suppose $M$ is toroidal. Since $M(\infty)$ does
not contain a non-separating torus or sphere, any essential torus
in $M$ must be separating. Let $F$ be an innermost essential torus
in $M$, i.e. $M=A\cup_F B$ where $\partial M \subseteq B$ and $A$
is atoroidal. Then $A, B$ are irreducible and
$\partial$-irreducible. By Lemma \ref{I2}, $F$ must be
compressible in $M(0)$, $M(\infty)$ and $M(1)$ and thus in $B(0)$,
$B(\infty)$ and $B(1)$. Applying Lemma \ref{generallemma} to $B$
with the two slopes $0 ,\infty$, we have two cases to consider.

\textit{Case 1: $B$ is a cable space $C(s,t)$, $t\geq 2$ with
cabling slope 0 or $\infty$}. Since $M(\infty)=S^2(2,3,5)$, which
does not have a lens space summand, $\infty$ cannot be the cabling
slope of $B$. Therefore $0$ is the cabling slope. Then
$B(0)=S^1\times D^2\# L(t,s)$ and $B(\infty)=S^1\times D^2$. Also
$B(1)=S^1\times D^2$ since $\Delta(0,1)=1$. Note that
$M(0)=A\cup_F B(0)$, $M(\infty)=A\cup_F B(\infty)$ and
$M(1)=A\cup_F B(1)$. If we let $r_0, r_\infty$ and $r_1$ be the
slopes on $F$ corresponding to the meridians of the solid tori of
$B(0)$, $B(\infty)$ and $B(1)$ respectively, then $M(0)=A(r_0)\#
L(t,s)$, $M(\infty)=A(r_\infty)$ and $M(1)=A(r_1)$. Lemma \ref{I2}
implies that $A(r_0)\cong L(3,1)$ or $L(4,1)$, $A(r_\infty)\cong
S^2(2,3,5)$, $A(r_1)\cong S^2(2,3,7)$ and $t=3$ or $4$.
Furthermore, from the proof of Lemma \ref{generallemma}, we know
that $\Delta(r_0, r_\infty)\geq t$. Thus $\Delta(r_0,
r_\infty)\geq 3$.

We have shown that the irreducible, $\partial$-irreducible and
atoroidal manifold $A$ admits two Dehn fillings $A(r_0)\cong
L(3,1)$ or $L(4,1)$, $A(r_\infty)\cong S^2(2,3,5)$ and
$\Delta(r_0, r_\infty)\geq 3$. This is a contradiction to
\cite[Theorem 1.1]{BZ1}, provided that $A$ is not Seifert fibered
(thus $A$ is hyperbolic).

To complete Case 1, we need to show that $A$ is not Seifert
fibered. Suppose $A$ is Seifert fibered. Then since
$A(r_\infty)\cong S^2(2,3,5)$ and $A(r_1)\cong S^2(2,3,7)$, $A$
must be homeomorphic to $D^2(2,3)$, which is homeomorphic to the
complement of the trefoil knot in $S^3$ i.e. the (2,3)-torus knot.
Also $A(\cong D^2(2,3))$ admits the lens space Dehn filling
$A(r_0)(=L(3,1)$ or $L(4,1))$. However this is impossible by \cite
[Lemma 7.4] {G2}.

\textit{Case 2: $B(0)$ and $B(\infty)$ are $S^1\times D^2$}. Note
that by applying Lemma \ref{generallemma} to two slopes $0$ and
$1$, and by knowing the fact that $M(1)$ doesn't have a lens space
summand, we can assume $B(1)$ to be $S^1\times D^2$. Let $r_0$,
$r_\infty$ and $r_1$ be the slopes on $F$ corresponding to the
meridians of $B(0)$, $B(\infty)$ and $B(1)$ respectively. Then
$M(0)=A(r_0)$ $(=L(3,1)\# L(4,1))$, $M(\infty)=A(r_\infty)$
$(=S^2(2,3,5))$ and $M(1)=A(r_1)$ $(=S^2(2,3,7))$. Also Lemma
\ref{generallemma} shows that $\Delta(r_0, r_\infty)\geq 4$. Then
it is easy to see by applying the same argument of the second
paragraph of this theorem that $A$ is not Seifert fibered and thus
hyperbolic. Then considering the manifold $A$, we get a
contradiction to \cite[Theorem 1]{BGZ}. This completes the proof.

\end{proof}

\end{document}